\documentclass{article}[12pt]
\usepackage{latexsym, amsmath, amsfonts, amssymb}
\usepackage{hyperref}

\newcommand{\nc}{\newcommand}
\nc{\ga}{\gamma} \nc{\di}{\displaystyle}
\nc{\ek}{\protect\\[1ex]}
\nc{\N}{{\mathbb N}} \nc{\R}{{\mathbb R}} \nc{\Z}{{\mathbb Z}}
\nc{\La}{\Lambda} \nc{\la}{\lambda} \nc{\da}{\delta}
\nc{\Da}{\Delta} \nc{\na}{\nabla} \nc{\vp}{\varphi} \nc{\si}{\sigma}
\nc{\Si}{\Sigma} \nc{\al}{\alpha} \nc{\be}{\beta} \nc{\om}{\omega}
\nc{\Om}{\Omega} \nc{\pa}{\partial} \nc{\ti}{\times}
 \nc{\ve}{\varepsilon} \nc{\ra}{\rightarrow} \nc{\Ra}{\Rightarrow}
\nc{\ran}{\rangle} \nc{\lan}{\langle}

 \nc{\eq}[1]{\mbox{\rm{(\ref{E#1})}}}
\nc{\qed}{\mbox{}\nolinebreak\hfill \rule{2mm}{2mm}}
 \nc{\ha}{\frac{1}{2}}
\nc{\hra}{\hookrightarrow} \nc{\supp}{{\rm supp}\,}
\nc{\curl}{\text{curl}\,} \nc{\dense}{\hra^{\hspace{-3mm}d\,\,}}

\newtheorem{lem}{Lemma}[section]
\newtheorem{theo}[lem]{Theorem}

\newtheorem{rem}[lem]{Remark}

\renewcommand{\div}{{\rm{div}}\,}

\numberwithin{equation}{section} 

\title{Energy Superposition and Regularity  for 3D Navier-Stokes
 Equations  in the Largest Critical Space}

\author{Myong-Hwan Ri\\
\small Institute of Mathematics, State Academy of Sciences, DPR
 Korea}
\date{}

  \allowdisplaybreaks[4]

\begin{document}
\bibliographystyle{alpha}

\maketitle%

\begin{abstract}
We show that  a Leray-Hopf weak solution to the 3D Navier-Stokes
Cauchy problem belonging to the space $L^\infty(0,T;
B^{-1}_{\infty,\infty}(\mathbb R^3))$ is regular in $(0,T]$. As a
consequence, it follows that any Leray-Hopf weak solution to the 3D
Navier-Stokes equations is regular while it is temporally bounded in
the largest critical space $\dot{B}^{-1}_{\infty,\infty}(\mathbb
R^3)$.

For the proof we present a new elementary method which is to
superpose the energy norm of high frequency parts in an appropriate
way to generate higher order norms. Thus, starting from the energy
estimates of high frequency parts of a weak solution, one can obtain
its estimates of higher order norms.

 By a linear energy superposition we get very simple and short proofs for
known regularity criteria for Leray-Hopf weak solutions in endpoint
Besov spaces $B^{\sigma}_{\infty,\infty}$ for $\sigma\in [-1,0)$,
the extension of Prodi-Serrin conditions. The main result of the
paper is proved by applying technique of a nonlinear energy
superposition and linear energy superpositions, repeatedly. The
energy superposition method developed in the paper can also be
applied to other supercritical nonlinear PDEs.

\end{abstract}

\noindent {\bf Keywords: } Navier-Stokes equations; regularity
criteria; largest critical space; energy superposition\\
{\bf 2010 MSC: } 35Q30, 76D05

\let\thefootnote\relax\footnote{\hspace{-0.3cm}
E-mail address\,$:$ math.inst@star-co.net.kp (Myong-Hwan Ri)}


\section{Introduction and main result}
Let us consider the Cauchy problem for the 3D Navier-Stokes
equations
\begin{equation}
\label{E1.1}
\begin{array}{rl}
     u_t-\nu\Da u  + (u\cdot\na)u+\na p  = 0 \,\, &\text{in }(0,\infty)\ti\R^3,\ek
      \div u = 0 \,\, &\text{in }(0,\infty)\ti\R^3,\ek
u(0,x)=u_0\, &\text{in }\R^3,
\end{array}
\end{equation}
where $u$, $p$ are velocity and pressure, respectively, and $\nu$ is
the kinematic viscosity of the observed viscous incompressible
fluid.

Leray proved in \cite{Le34} that if $u_0\in L^2(\R^3)$, $\div
u_0=0$, then the problem \eq{1.1} has a global weak solution
(Leray-Hopf weak solution)
$$
\label{E1.3} u\in L^2(0,\infty; H^1(\R^3))\cap L^\infty(0,\infty;
L^2(\R^3)),
 $$
  which satisfies \eq{1.1} in a weak sense and the {\it energy inequality}
\begin{equation}
\label{EEI}
\begin{array}{l}\di\frac{1}{2}\|u(t)\|_2^2+\nu\int_0^t\|\na
u(\tau)\|_2^2\,d\tau\leq \frac{1}{2}\|u_0\|_2^2,\; \forall t\in
(0,\infty).
 \end{array}
\end{equation}
 A weak solution $u$  to \eq{1.1} is called regular (or
strong) in  $(0,T]$, $0<T<\infty$, if
\begin{equation}
\label{E1.4n}
 u\in L^\infty(0,T; H^1(\R^3))\cap L^2(0,T; H^2(\R^3)).
  \end{equation}

The problem on global uniqueness and regularity of Leray-Hopf weak
solutions has a very long history. In their famous works
\cite{Pr59}, \cite{Se63}, Prodi and Serrin obtained the condition
 \begin{equation}
 \label{E1.4}
u\in L^p(0,T;L^q(\R^3)), \frac{2}{p}+\frac{3}{q}=1, 2\leq p<\infty,
3<q \leq \infty,
 \end{equation}
 for a Leray-Hopf weak solution $u$ to be regular in $(0,T]$.
 The sufficient regularity condition \eq{1.4} was extended to the endpoint Besov space
 version
 \begin{equation}
 \label{E1.5}
u\in L^p(0,T; \dot{B}^\si_{\infty,\infty}(\R^3)), \frac{2}{p}=1+\si,
 \end{equation}
for the case $2<p<\infty$ by Kozono and Shimada \cite{KoSh04} (see
Cheskidov and Shvydkoy \cite{CS10} where the homogeneous Besov space
$\dot{B}^\si_{\infty,\infty}(\R^3)$ is extended to the inhomogeneous
Besov space $B^\si_{\infty,\infty}(\R^3)$), and for the case
$1<p\leq 2$, for example, by Chen and Zhang \cite{ChZh06}. See
\cite{ZY16} for vorticity conditions corresponding to \eq{1.5},
i.e., $\curl u\in L^p(0,T; \dot{B}^{\si-1}_{\infty,\infty}(\R^3))$.

 For the case $p=1$, Beale, Kato and Majda \cite{BKM84} obtained  regularity condition
 $$\curl u \in L^1(0,T;L^\infty(\R^3)),$$
 which was extended to the condition $\curl u \in L^1(0,T; BMO)$ in
\cite{KoTa00} using the logarithmic Sobolev inequality in Besov
spaces, and to the still weaker condition $\curl u\in
L^1(0,T;\dot{B}^0_{\infty,\infty}(\R^3))$ (\cite{KoOgTa02}).

\par\medskip

Extending to the limiting case $p=\infty$ in \eq{1.4} is generally
regarded as difficult. By exploiting the backward uniqueness
property of parabolic equations, Eskauraiza, Seregin and Sv\'erak in
\cite{ESS03}(2003) showed that a Leray-Hopf weak solution is regular
in $(0,T]$ if
 \begin{equation}
 \label{E1.6}
u\in L^\infty(0,T; L^3(\R^3)).
 \end{equation}
 Alternatively,
 Gallagher, Koch and Planchon in \cite{GKB13}(2013) proved for
  a mild solution $u$ to \eq{1.1} with initial values in $L^3(\R^3)$ that a potential
 singularity at $t=T$ implies $\lim_{t\ra
 T-0}\|u(t)\|_3=+\infty$, and  in
\cite{GKB16}(2016) extended the result for $L^3(\R^3)$ to wider
critical Besov spaces $\dot{B}^{-1+3/p}_{p,q}(\R^3)$, $3<p,q<\infty$
 by developing a profile decomposition technique and using the method of
``critical elements" developed in \cite{KeMe06}-\cite{KeMe10}
 (cf. \cite{ChPl13}(2013) for the case $3<p<\infty, q<2p/(p-1)$). We also
refer the readers to \cite{Ph15}(2015), where the result of
\cite{ESS03} was directly extended to the Lorentz space
$L^{3,q}(\R^3)$, $3<q<\infty$; see \cite{WaZh17}(2017) for further
improvement towards the space $BMO^{-1}$. Whether boundedness of a
NS weak solution in the largest (i.e., the weakest) critical space
$\dot{B}^{-1}_{\infty,\infty}(\R^3)$ will imply its regularity is a
very important question ({\it say},
 $L^\infty(\dot{B}^{-1}_{\infty,\infty})$-regularity problem) for analyzing the
 role of critical norms in the regularity theory; if the problem be solved positively,
 then the problem on regularity conditions for whole velocity
 in critical spaces will be ultimately settled.
 We recall that a
critical space for the 3D Navier-Stokes equations is the space whose
norm is invariant by the scaling $u_\la(x)=\la u(\la x)$, $\la>0$.
The following embedding holds
 between critical spaces for the 3D Navier-Stokes equations:
 \begin{equation}
 \label{E1.8}
\begin{array}{l}
   L^3\hra \dot{B}^{-1+3/p}_{p,\infty}\hra
  BMO^{-1}
  \hra \dot{B}^{-1}_{\infty,r}\hra
  \dot{B}^{-1}_{\infty,\infty},\; 3 \leq p<\infty,2<r<\infty.
  \end{array}
  \end{equation}
It is remarkable that \eq{1.1} is well-posed for small initial
values in the critical spaces $L^3(\R^3)$ (\cite{Ka84}) and
$BMO^{-1}$ (\cite{KoTa01}), whereas the problem is ill-posed in
$\dot{B}^{-1}_{\infty,\infty}(\R^3)$ (\cite{BouPa08}) and even in
$\dot{B}^{-1}_{\infty,q}(\R^3)$, $q\in [1,\infty]$ (\cite{Wa15}) in
the sense that a norm inflation occurs; note that
$\dot{B}^{-1}_{\infty,q}(\R^3)\hra BMO^{-1}$ for $1\leq q \leq 2$.

Concerning the above-mentioned
$L^\infty(\dot{B}^{-1}_{\infty,\infty})$-regularity problem,
Cheskidov and Shvydkoy \cite{CS10}(2010), using Paley-Littlewood
theory, showed that a Leray-Hopf weak solution $u$ to \eq{1.1} is
regular in $(0,T]$ if
 $$\lim_{k\ra\infty}\|u^{k}\|_{L^\infty(0,T;
B^{-1}_{\infty,\infty}({\mathbb R}^3))}<C\nu, $$
 or if $u(T)\in
B^{-1}_{\infty,\infty}(\R^3)$, $\lim_{t\ra
T-0}\|u(t)-u(T)\|_{B^{-1}_{\infty,\infty}}<C\nu$ with some absolute
constants $C>0$,  where and in what follows we use the notation
\begin{equation}
\label{E2.1}
\begin{array}{l}
 u^k(t,x):=\frac{1}{(2\pi)^{3/2}}\int_{|\xi|\geq k}{\mathcal F}{u}(t,
\xi)e^{ix\cdot\xi}\,d\xi,\ek u_k:=u-u^k,\;
  u_{h,k}:=u^{h}-u^{k}\quad \text{for } 0\leq h<k,
  \end{array}
\end{equation}
and ${\mathcal F}{u}(\xi)\equiv
\widehat{u}(\xi):=\frac{1}{(2\pi)^{3/2}}\int_{\R_\xi^3}u(x)e^{-ix\cdot\xi}dx$
stands for the three-dimensional Fourier transform of $u$. In
\cite{Ri19}(2019) and \cite{Ri22}(2022), the author showed that if
$L^\infty(0,T; {B}_{\infty,\infty}^{-1})$-norm or $L^1(0,T;
{B}_{\infty,\infty}^{-1})$- or $L^2(0,T;
{B}_{\infty,\infty}^{-1})$-norm of a mid frequency part $u_{k/2,k}$
is small depending on $\nu$, $u_0$, $k$ and
 $$A:=\sup_{t\in (0,T)} \big(\nu t^{-1}\int_0^t \|\na u\|_2^2d\tau\big)$$
 for some $k\geq k_*(\nu,u_0,A)$, then regularity of $u$ in $(0,T]$ follows.

\par\medskip
In this paper, we give the positive answer to the
$L^\infty(\dot{B}^{-1}_{\infty,\infty})$-regularity problem; we even
 show that $L^3(\R^3)$ in \eq{1.6} can be
extended to the space $B^{-1}_{\infty,\infty}(\R^3)$ which is larger
than $\dot{B}^{-1}_{\infty,\infty}(\R^3)$. The main result of the
paper is as follows.
\begin{theo}
\label{T1.1} {\rm Let $u$ be a Leray-Hopf weak solution to \eq{1.1}
with $u_0\in H^{1}(\R^3)$, $\div u_0=0$. If
 \begin{equation}
 \label{E1.9}
u\in L^\infty(0,T; {B}^{-1}_{\infty,\infty}(\R^3)),
 \end{equation}
then
$$u\in
L^\infty(0,T;H^{1}(\R^3))\cap L^2(0,T;H^{2}(\R^3))
$$
and hence $u$ is regular in $(0,T]$.
 }
\end{theo}
 \begin{rem}
The condition \eq{1.9} may be the last extended
 Serrin-type  regularity condition for velocity itself.
Moreover, Theorem \ref{T1.1} has the following important
implications:

\begin{itemize}
 \item[(i)] Boundedness in any critical norm of \eq{1.8} for a Leray-Hopf weak
solution to \eq{1.1} implies its regularity.

\item[(ii)] If a Leray-Hopf weak solution $u$ is irregular at time $T$,
then its weakest critical norm
$\|u(t)\|_{\dot{B}^{-1}_{\infty,\infty}(\R^3)}$ should be unbounded
near the epoch $T$. In other words,  there should exist a sequence
of wavenumbers $\{k_n\}_{n\in\N}$ such that as $k_n\nearrow\infty$
the velocity part of around the wavenumber $k_n$ should grow faster
than any linear growth in $k_n$. It means that, in the possible
blow-up scenario, the convection term $(u\cdot\na)u$ should always
dominate over the viscosity term $-\Da u$ since $(u\cdot\na)u\sim
kU^2$ and $-\Da u\sim k^2U$, where $U$ is the characteristic
velocity corresponding to the wavenumber $k$.
 \end{itemize}
 \end{rem}

In order to prove the main result we present a new method, which is
based on energy estimates of  high frequency parts of the weak
solution and their appropriate superposition.

The energy method is successful for showing the global existence of
weak solutions to Navier-Stokes equations thanks to the crucial
cancelation property
 $((v\cdot\na)u,u)=0$  for solenoidal vector field $v$.
  However, the method has been regarded as
ineffective for global existence of the strong solution and global
regularity of weak solutions, since the conventional arguments are
usually based on testing the equations with $(-\Da)^{\al}u$ for
$\al\geq 1/2$ and then the estimate of $((u\cdot\na)u,(-\Da)^{\al}
u)$ involves essential difficulty; the situation does not change
when applying the energy method to the vorticity equation
$$\om_t-\nu\Da \om  + (u\cdot\na)\om-(\om\cdot\na)u  = 0$$
for $\om=\curl u$ since $((\om\cdot\na)u,\om)=((u\cdot\na)u,\Da u)$.

When the energy method is applied to high frequency parts of
\eq{1.1}, the nonlinear term is not completely canceled but the
partial (but subtle) cancelation $((u\cdot\na)u^k,u^k)=0$  makes it
possible to obtain the estimate
 \begin{equation}
 \label{E1.11}
 \frac{d}{dt}\|u^k\|_2^2+\nu \|\na
u^k\|_2^2\lesssim
k^{1-\si}\|u_k\|_{B^{\si}_{\infty,\infty}}\|u^{k/2}\|_2^2,\quad
\si\in [-1,0),k\geq 1.
 \end{equation}
  This
differential inequality for the energy norm of high frequency parts
 is the starting point of our method.
In order to obtain higher order norm estimate of the weak solution
starting from its zeroth order estimate \eq{1.11}, we observe the
relation
 \begin{equation}
 \label{E1.7}
\|u\|_{\dot{H}^{1/2}}^2\sim \sum_{k\in\N\cup\{0\}}\|u^k\|_2^2,
 \end{equation}
see \eq{2.11} for more details.
  The right-hand side of \eq{1.7} may be considered, in a sense, to be
a {\it linear superposition} of energy of high frequency parts
(i.e., the sum of energy for high frequency parts with equally
increasing lower wavenumber bounds). By the direct use of \eq{1.11}
and \eq{1.7}, known regularity criteria for
 Leray-Hopf weak solutions in endpoint Besov spaces
$B^{\si}_{\infty,\infty}$ or $\dot{B}^{\si}_{\infty,\infty}$,
$\si\in [-1,0)$, the extension of Prodi-Serrin
 conditions can be obtained with short proofs,
against the  ever known long proofs with usually complicated
arguments in a number of articles (Section 2).

Theorem \ref{T1.1} is proved using more or less delicate techniques
of superposition, being based on the fact that any higher order norm
can be generated by a possible nonlinear superposition (i.e. the sum
of energy for high frequency parts with unequally increasing lower
wavenumber bounds) of energy of high frequency parts. Starting from
\eq{1.11}, we superpose energy of high frequency parts with
polynomially increasing lower wavenumber bound; by regulating the
polynomial order dependent upon the size of
$L^\infty(B^{-1}_{\infty,\infty})$-norm, we are led to a
differential inequality with a small coefficient of the major term
on the right-hand side. Then, we apply the technique of a linear
superposition to the resulting inequality finite times repeatedly,
obtaining a required higher order norm bound of the weak solution
(Section 3).

\begin{rem}
The key ingredients of the proposed method are the energy estimates
of high frequency parts and appropriate superposition of their
energy. Therefore, we would call the method ``{\em energy
superposition method}".
\end{rem}

We use the following notations. The sets of all natural numbers and
all integers are denoted by $\N$ and $\Z$, respectively. We denote
by $H^s(\R^3)$ and $\dot{H}^s(\R^3)$ for  $s>0$ the Bessel potential
and homogeneous Bessel potential spaces, respectively.  The Besov
spaces and homogeneous Besov spaces are denoted by $B_{p,q}^s(\R^3)$
and $\dot{B}_{p,q}^s(\R^3)$, $s\in\R, 1\leq p,q\leq \infty$,
respectively.  We do not distinguish between scalar and vectorial
function spaces, and denote the generic constants by the same letter
$c,C$ unless specially necessitated.

\section{Preliminaries}

Let
$$\Da_j u= (\chi_j(\xi)\hat{u}(\xi))^\vee\;\text{for }j\in\Z,
\;\tilde\Da_0u=(\tilde\chi_0(\xi)\hat{u}(\xi))^\vee, u\in {\cal
S}'(\R^n),$$ where $\chi_j$, $j\in\Z$, and $\tilde\chi_0$ are the
characteristic functions of the sets $\{\xi\in\R^3:  2^{j-1}\leq
|\xi|< 2^{j}\}$ and $\{\xi\in\R^3: |\xi|< 1\}$, respectively, and
${\cal S}'(\R^3)$ is the space of tempered distributions. Then,
$\{\Da_j\}_{j\in\Z}$ becomes a dyadic Schauder decomposition for
homogeneous Besov spaces $\dot{B}^s_{q,r}(\R^3)$, $s\in\R, 1\leq
q,r\leq\infty$, i.e., $\Da_j\Da_l=0$ for $j\neq l$ and
$u=\sum_{j\in\Z} \Da_j u$ for all $u\in \dot{B}^s_{q,r}(\R^3)$ such
that the sum converges unconditionally in $\dot{B}^s_{q,r}(\R^3)$.
Meanwhile, $\{\tilde\Da_0,\Da_j\}_{j\in \N}$ becomes a dyadic
Schauder decomposition
 for inhomogeneous Besov spaces $B^s_{q,r}(\R^3)$, $s\in\R, 1\leq q,r\leq\infty$.

Recall that for $s\in\R, 1\leq q\leq \infty$, $1\leq r<\infty$ the
norms
\begin{equation}
 \label{E3.1}
(\sum_{j\in\Z}\|2^{sj}\Da_j u\|^r_{L^q(\R^3)})^{1/r}, \quad
\sup_{j\in\Z}\|2^{sj}\Da_j u\|_{L^q(\R^3)}
 \end{equation}
  are equivalent to the
homogeneous Besov norms $\|u\|_{\dot{B}_{q,r}^s(\R^3)}$,
$\|u\|_{\dot{B}_{q,\infty}^s(\R^3)}$, respectively, introduced via
Paley-Littlewood decomposition. In the same context,
  \begin{equation}
 \label{E3.2}
 (\|\tilde\Da_0u\|^r_{L^q(\R^n)}+\sum_{j\in\N}\|2^{sj}\Da_j u\|^r_{L^q(\R^3)})^{1/r},
  \;
  \max\{\|\tilde{\Da}_0
u\|_{L^q(\R^3)}, \sup_{j\in\N}\|2^{sj}\Da_j u\|_{L^q(\R^3)}\}
 \end{equation}
 are equivalent to the inhomogeneous Besov norms
 $\|u\|_{B_{q,r}^s(\R^3)}$,
$\|u\|_{B_{q,\infty}^s(\R^3)}$, respectively, introduced via
Paley-Littlewood decomposition. Henceforth, in the sequel, we
suppose that the homogeneous and inhomogeneous Besov spaces are
endowed with the norms defined by \eq{3.1} and \eq{3.2},
respectively.

\medskip
In the below, let $u\in L^2(\R^n)$ be suitably smooth so that the
appearing norms in lines make senses. Following the notation
\eq{2.1}, one has
$$\Da_j u=u_{2^{j-1},2^j}, j\in\Z,\quad \tilde{\Da}_0 u=u_1.$$

It is easy, by Plancherel's equality, to see that
 \begin{equation}
 \label{E3.3}
(u_k, u^l)_2=0,\quad 0<k\leq l<\infty
 \end{equation}
  and
  \begin{equation}
 \label{E3.4}
k^\al\|u^k\|_2^2\leq \|u^k\|_{\dot{H}^{\al}}^2,\;\forall k>0, \al>0.
 \end{equation}
Obviously,
  \begin{equation}
 \label{E3.5}
\supp{\widehat{(u_ku_l)}}\subset \{\xi\in\R^3: |\xi|\leq
k+l\},\;\forall k,l>0.
 \end{equation}

It follows that
 \begin{equation}
 \label{E2.21}
  \|u_k\|_\infty\leq c(\si)k^{-\si}\|u_k\|_{B^\si_{\infty,\infty}},\quad \forall \si\in [-1,0),\forall k\geq 1.
\end{equation}
In fact, for $\si\in [-1,0)$ we have
 $$\begin{array}{rl}
 \|u_k\|_\infty&\leq \|u_1\|_\infty+
 \di\sum_{j=1}^{j_0}\|[u_k]_{2^{j-1},2^{j}}\|_\infty\ek
 &\leq \|u_1\|_\infty+
 c\di\sum_{j=1}^{j_0} 2^{-j\si}\big\||\na|^\si
 [u_k]_{2^{j-1},2^{j}}\big\|_\infty\ek
&\leq
 k^{-\si}\|u_1\|_\infty+c\di\big(\sum_{j=1}^{j_0} 2^{-j\si}\big)
   \big\||\na|^\si u_k\big\|_{B^{0}_{\infty,\infty}}\ek
&\leq k^{-\si}\|u_1\|_\infty+ c 2^{-j_0\si+1}
   \big\||\na|^\si u_k\big\|_{B^{0}_{\infty,\infty}}
\ek&\leq
 c(\si) k^{-\si}\big\|u_k\big\|_{B^{\si}_{\infty,\infty}},
 \end{array}$$
 where $|\na|:=(-\Da)^{1/2}$ and $j_0\in\Z$ is such that
  $2^{j_0-1}\leq k < 2^{j_0}$.
 Thus, the inequality \eq{2.21} holds true.

  Furthermore, we have
  \begin{equation}
 \label{E3.6}
\|u_{k,l}\|_{\dot{B}^{-1}_{\infty,\infty}}\leq
\|u^{k}\|_{\dot{B}^{-1}_{\infty,\infty}},\quad \forall l>k>0.
 \end{equation}
In fact, in view of the definition \eq{3.1} for homogeneous Besov
norms, \eq{3.6} is trivial when $l=2^j$ for some $j\in\Z$. Let $l\in
(2^{j-1},2^j)$. For the case $k\geq 2^{j-1}$ we have
 \begin{equation}
 \label{E3.7}
\begin{array}{l}
\|u_{k,l}\|_{\dot{B}^{-1}_{\infty,\infty}}
=2^{-j}(2\pi)^{-\frac{3}{2}}\|u_{k,l}\|_\infty\\[3ex]\quad\di
 =2^{-j}(2\pi)^{-\frac{3}{2}}\big\|\int_{k<|\xi|<l}
\widehat{u_{k,l}}(\xi)e^{ix\cdot\xi}d\xi\big\|_{L_x^\infty(\R^3)}\\[3ex]
 \quad \di
 = 2^{-j}(2\pi)^{-\frac{3}{2}}(2^{-j}l)\Big\|\int_{(2^{j}k)/l<|\eta|<2^{j}}
\widehat{u_{k,l}}\big((2^{-j}l)^{1/3}\eta\big)
      e^{ix\cdot (2^{-j}l)^{1/3}\eta} d\eta\Big\|_{L_x^\infty(\R^3)}\\[3ex]
 \quad \di
 \leq 2^{-j}(2\pi)^{-\frac{3}{2}}(2^{-j}l)\Big\|\int_{(2^{j}k)/l<|\eta|<2^{j}}
 \widehat{v}(\eta) e^{i(2^{-j}l)^{1/3}x\cdot\eta}d\eta\Big\|_{L_x^\infty(\R^3)}\\[3ex]
 \quad\di = 2^{-j} (2^{-j}l)\Big\|
 v\big((2^{-j}l)^{1/3}\cdot\big)\Big\|_\infty
  \; = (2^{-j}l)\Big\|v\big((2^{-j}l)^{1/3}\cdot\big)\Big\|_{\dot{B}^{-1}_{\infty,\infty}}\\[3ex]
  \quad\di = (2^{-j}l)^{2/3}\big\|v\big\|_{\dot{B}^{-1}_{\infty,\infty}}
           \leq
           (2^{-j}l)^{2/3}\big\|V\big\|_{\dot{B}^{-1}_{\infty,\infty}},
 \end{array}
 \end{equation}
where
$\widehat{v}(\eta):=\widehat{u_{k,l}}\big((2^{-j}l)^{1/3}\eta\big)$,
 $\widehat{V}(\eta):=\widehat{u^{k}}\big((2^{-j}l)^{1/3}\eta\big)$
 and we used the critical property of the $\dot{B}^{-1}_{\infty,\infty}$-norm.
By the reverse argument we have
 \begin{equation}
 \label{E3.8}
\begin{array}{l}\big\|V\big\|_{\dot{B}^{-1}_{\infty,\infty}}=
 \|{\cal F}^{-1}\widehat{u^k}\big((2^{-j}l)^{1/3}\eta\big)\|_{\dot{B}^{-1}_{\infty,\infty}}
 \ek\qquad =(2^{-j}l)^{-1}\|u^k\big((2^{-j}l)^{-1/3}\cdot\big)\|_{\dot{B}^{-1}_{\infty,\infty}}$$
 \ek\qquad =(2^{-j}l)^{-2/3}\|u^k\|_{\dot{B}^{-1}_{\infty,\infty}}.$$
  \end{array}
  \end{equation}
Consequently, \eq{3.7} and \eq{3.8} imply \eq{3.6}. For the last
case $k<2^{j-1}$ we have, by above proved fact,
$$\begin{array}{l}\|u_{k,l}\|_{\dot{B}^{-1}_{\infty,\infty}}=
 \max\{\|u_{k,2^{j-1}}\|_{\dot{B}^{-1}_{\infty,\infty}},
  \|u_{2^{j-1},l}\|_{\dot{B}^{-1}_{\infty,\infty}}\}\ek
\qquad \leq \max\{\|u_{k,2^{j-1}}\|_{\dot{B}^{-1}_{\infty,\infty}},
  \|u^{2^{j-1}}\|_{\dot{B}^{-1}_{\infty,\infty}}\}=\|u^{k}\|_{\dot{B}^{-1}_{\infty,\infty}}.
   \end{array}$$
Thus, \eq{3.6} is proved.

Note that
$\|u^k\|_{\dot{B}^{-1}_{\infty,\infty}}=\|u^k\|_{B^{-1}_{\infty,\infty}}$
for $k\geq 1$. Hence, it follows from \eq{3.6} that
  \begin{equation}
 \label{E3.9}
\|u_{k,l}\|_{B^{-1}_{\infty,\infty}}\leq
\|u^{k}\|_{B^{-1}_{\infty,\infty}},\quad \forall l>k>0.
 \end{equation}

\section{Short proofs for regularity criteria in endpoint Besov spaces }

In this section, introducing the {\it energy superposition method},
we give short proofs for already known regularity criteria for 3D
Navier-Stokes equations.

\begin{theo}
 \label{T2.2}
 {\rm Let $u$ be a Leray-Hopf weak solution  to \eq{1.1} with $u_0\in H^1(\R^3)$, $\div
 u_0=0$.
 If
 $$u\in L^{2/(1+\si)}(0,T;B^\si_{\infty,\infty}(\R^3))$$
   for some $\si\in (-1,0)$, then $u$ is regular in $(0,T]$.
 }
\end{theo}

{\bf Proof:} Let us assume without loss of generality that $T$ is
the first blow-up epoch for $u$. By $L^2$-scalar product to \eq{1.1}
with $u^k$, $k>0$, we have
 \begin{equation}
 \label{E2.18}
\begin{array}{l}
\frac{d}{2dt}\|u^k\|_2^2+\nu\|\na u^k\|_2^2=-((u\cdot\na)u_k,
u^k)\ek
 \quad =-((u_k\cdot\na)u_k-(u^k\cdot\na)u_k, u^k)\ek
 \quad
=-((u_k\cdot\na)u_{k/2,k},u^{k})-((u_{k/2,k}\cdot\na)u_{k/2},u^{k})-(u^k\cdot\na)u_k,u^k)
 \;\text{in } (0,T)
\end{array}
 \end{equation}
  in view of  $u_k=u_{k/2}+u_{k/2,k}$ and
\eq{3.5}. In the right-hand side of \eq{2.18} we have
$$\begin{array}{rcl}
|((u_{k/2,k}\cdot\na)u_{k/2},u^{k})+(u^k\cdot\na)u_k,u^k)|
 &\leq & (\|\na u_{k/2}\|_\infty+\|\na u_k\|_\infty)\|u^{k/2}\|_2^2\ek
 &\leq & ck\|u_k\|_\infty\|u^{k/2}\|_2^2
  \end{array}$$
  by H\"older's inequality.
 On the other hand, we have
$$ \begin{array}{rcl} |((u_k\cdot\na)u_{k/2,k},u^{k})|
 &\leq &  c\|u_k\|_\infty\|\na u_{k/2,k}\|_2\|u^k\|_2\ek
  &\leq & c k\|u_k\|_\infty\|u^{k/2}\|_2^2.
 \end{array}
 $$
 Therefore, we get from \eq{2.18} that
 \begin{equation}
 \label{E2.19}
 \begin{array}{l}
 \frac{d}{2dt}\|u^k\|_2^2+\nu\|\na u^k\|_2^2
      \leq ck\|u_k\|_\infty\|u^{k/2}\|_2^2
      \quad\text{in }(0,T).
      \end{array}
 \end{equation}

 Using  \eq{2.21} and \eq{3.4}, we get from \eq{2.19} that, ,
 \begin{equation}
 \label{E2.15}
  \begin{array}{rcl}
 \frac{d}{2dt}\|u^k\|_2^2+\nu\|\na u^k\|_2^2
      &\leq & c k^{1-\si}\|u_k\|_{B^{\si}_{\infty,\infty}}\|u^{k/2}\|_2^2\ek
      & \leq & c \|u_k\|_{B^{\si}_{\infty,\infty}}\||\na|^{(1-\si)/2}u^{k/2}\|_2^2\ek
      & \leq & c \|u_k\|_{B^{\si}_{\infty,\infty}}\|u^{k/2}\|_2^{1+\si}\|\na
      u^{k/2}\|_2^{1-\si}\ek
   & \leq & \frac{\nu}{2}\|\na
      u^{k/2}\|_2^{2}+c(\nu,\si)
      \|u_k\|_{B^{\si}_{\infty,\infty}}^{2/(1+\si)}\|u^{k/2}\|_2^{2}
      \end{array}
 \end{equation}
 for $\si\in (-1,0)$,
where an interpolation inequality is used.
 \par\medskip

Adding up \eq{2.15} over $k\in\N$ yields
$$\begin{array}{l} \di
\frac{d}{2dt}\sum_{k\in\N}\|u^{k}\|_2^2+\nu \sum_{k\in\N}\|\na
u^k\|_2^2\leq \frac{\nu}{2} \sum_{k\in\N}\|\na u^{k/2}\|_2^2
+c\sum_{k\in\N}\|u^{k/2}\|_2^2\|u_k\|^{2/(1+\si)}_{B^{\si}_{\infty,\infty}}\\[3ex]
      \leq \di\nu \sum_{k\in\N}\|\na u^{k}\|_2^2+\frac{\nu}{2}\|\na u^{1/2}\|_2^2
+c\big(2\|u\|^{2/(1+\si)}_{B^{\si}_{\infty,\infty}}
 \sum_{k\in\N}\|u^{k}\|_2^2+\|u^{1/2}\|_2^2\|u_1\|^{2/(1+\si)}_{B^{\si}_{\infty,\infty}}\big),
\end{array}$$
where the constants $c$ depend only on $\nu$ and $\si$. Hence, using
that $\|u_1\|_{B^{\si}_{\infty,\infty}}\leq c\|u_1\|_2$ and the
energy inequality, we are led to
\begin{equation}
\label{E2.2} \frac{d}{dt}\sum_{k\in\N}\|u^{k}\|_2^2\leq
c\|u\|^{2/(1+\si)}_{B^{\si}_{\infty,\infty}}\sum_{k\in\N}\|u^{k}\|_2^2
 +\frac{\nu}{2}\|\na u^{1/2}\|_2^2+e(\si,\nu,\|u_0\|_2).
 \end{equation}
Hence,
$$\begin{array}{l}
\di\sum_{k\in\N}\|u^{k}(t)\|_2^2
 \leq
c\int_0^t\|u(\tau)\|^{2/(1+\si)}_{B^{\si}_{\infty,\infty}}
 \sum_{k\in\N}\|u^{k}(\tau)\|_2^2\,d\tau\ek\qquad\quad
 +\di\sum_{k\in\N}\|u_0^{k}\|_2^2+\tilde{e}(\si,\nu, \|u_0\|_2,T),\quad\forall t\in
 (0,T),
 \end{array}$$
 which yields by Gronwall's inequality that
\begin{equation}
 \label{E2.10}
\sum_{k\in\N}\|u^{k}(t)\|_2^2\leq \big(\tilde{e}(\si,\nu,
\|u_0\|_2,T)+\sum_{k\in\N}\|u_0^{k}\|_2^2\big)
  \exp\{c\|u\|^{2/(1+\si)}_{L^{2/(1+\si)}(0,T;B^{\si}_{\infty,\infty})}\}
  \end{equation}
  for all $ t\in (0,T)$.
 Note that
 \begin{equation}
 \label{E2.11}
 \begin{array}{l}
 \di\sum_{j=k}^\infty\| u^{j}\|_2^2=\sum_{j=k}^\infty(j-k+1)\| u_{j,j+1}\|_2^2
  =\sum_{j=k}^\infty j\|u_{j,j+1}\|_2^2-(k-1)\|u^k\|_2^2,\ek
 \di\sum_{j=k}^\infty j\|u_{j,j+1}\|_2^2\sim
\|u^k\|_{\dot{H}^{1/2}}^2,\quad \forall k\in\N.\end{array}
\end{equation}
Therefore, \eq{2.10} and \eq{2.11} with $k=1$ yield that
$$\|u^{1}(t)\|_{\dot{H}^{1/2}}^2\leq \tilde{e_1}(\si,\nu,
\|u_0\|_{H^{1/2}},T)
\exp\{c\|u\|^{2/(1+\si)}_{L^{2/(1+\si)}(0,T;B^{\si}_{\infty,\infty})}\}.
$$
Thus, we have $u\in L^\infty(0,T; \dot{H}^{1/2}(\R^3))$, and $u$ can
not blow up at $t=T$. The proof is complete. \hfill\qed

\begin{rem}
When $\si$ moves in $(-1,0)$, the temporal integral exponent
$2/(1+\si)$ in the statement of Theorem \ref{T2.2} runs through
$(2,\infty)$. Hence, Theorem \ref{T2.2} includes the results of
\cite{KoSh04}, \cite{CS10}. For the case $\si=-1$, in the below, we
also have much shorter and simpler proof for the same result than
\cite{CS10}.
\end{rem}

\begin{theo}
 \label{T2.3}
 {\rm Let $u$ be a Leray-Hopf weak solution $u$ to \eq{1.1} with $u_0\in H^1(\R^3)$, $\div
 u_0=0$.
 If
  $$\|u^{k}\|_{L^\infty(0,T;{B}^{-1}_{\infty,\infty}(\R^3))}\leq c\nu$$
  holds with some $k>0$ and absolute constant $c>0$,
     then  $u$ is regular in $(0,T]$.
 }
 \end{theo}
 {\bf Proof:} Let us assume without loss of generality that $T$ is
the first blow-up epoch for $u$. Let us fix $k\in\N$.
 Using \eq{2.21} with $\si=-1$ and \eq{3.9}, we get from \eq{2.19}
 that
 \begin{equation}
\label{E2.12}
\begin{array}{l}
 \di\frac{d}{2dt}\|u^j\|_2^2+\nu\|\na u^j\|_2^2
    \leq
    cj(\|u_k\|_\infty+\|u_{k,j}\|_\infty)\|u^{j/2}\|_2^2,
\ek \qquad\leq
    c(j\|u_k\|_\infty+j^2\|u_{k,j}\|_{B^{-1}_{\infty,\infty}})\|u^{j/2}\|_2^2
    \ek
\qquad\leq
    c(j\|u_k\|_\infty+j^2\|u^{k}\|_{B^{-1}_{\infty,\infty}})\|u^{j/2}\|_2^2,
     \quad \forall j\geq k,
  \end{array}
\end{equation}
with some absolute constant $c>0$.
Then, by adding up \eq{2.12} over $j\geq k$, we have
 \begin{equation}
\label{E2.14}
\begin{array}{l}
 \di\frac{d}{2dt}\sum_{j\geq k}\| u^{j}\|_2^2+\nu\sum_{j\geq k}\|\na
u^{j}\|_2^2\ek \di
    \lesssim \|u_k\|_\infty \sum_{j\geq k}j\|u^{j/2}\|_2^2
              +\|u^k\|_{B_{\infty,\infty}^{-1}}\sum_{j\geq k}j^2\|u^{j/2}\|_2^2\ek
  \di\lesssim  \|u_k\|_\infty\Big(\sum_{k\leq j< 2k}j\|u^{j/2}\|_2^2
               +\sum_{j\geq 2k}j\|u^{j/2}\|_2^2\Big)
               \ek\qquad\qquad\di
          + \|u^k\|_{B_{\infty,\infty}^{-1}}\Big(\sum_{k\leq j<2k}j^2\|u^{j/2}\|_2^2+
           \sum_{j\geq 2k}\|\na u^{j/2}\|_2^2\Big)\ek
\di\lesssim
 k^{7/2}\|u_0\|_2^3+k^{3/2}\|u_0\|_2\sum_{j\geq
    k}j\|u^{j}\|_2^2
     +k^3\|u^k\|_{B_{\infty,\infty}^{-1}}\|u_0\|_2^2
         +\|u^k\|_{B_{\infty,\infty}^{-1}}\sum_{j\geq
    k}\|\na u^{j}\|_2^2
  \end{array}
\end{equation}
in view of
$$\|u_k\|_\infty\lesssim k^{3/2}\|u_0\|_2, \sum_{j\geq
    2k}j\|u^{j/2}\|_2^2\leq 2\sum_{j\geq
    k}j\|u^{j}\|_2^2, \sum_{j\geq
    2k}\|\na u^{j/2}\|_2^2\leq 2\sum_{j\geq
    k}\|\na u^{j}\|_2^2$$ and the energy
inequality. Note that, by \eq{3.4},
$$\sum_{j\geq
    k}j\|u^{j}\|_2^2\leq \ve \sum_{j\geq
    k}j^2\|u^{j}\|_2^2+c(\ve)\sum_{j\geq
    k}\|u^{j}\|_2^2\leq \ve \sum_{j\geq
    k}\|\na u^{j}\|_2^2+c(\ve)\sum_{j\geq
    k}\|u^{j}\|_2^2.$$
Therefore, when
$$\|u^k\|_{L^\infty(0,T;B_{\infty,\infty}^{-1})}<C\nu$$
for an absolute constant $C>0$, it follows from \eq{2.14} for
sufficiently small $\ve>0$ depending on $\nu,\|u_0\|_2$ and $k$ that
$$\di\frac{d}{dt}\sum_{j\geq k}\| u^{j}\|_2^2\leq c(k,\|u_0\|_2,\nu)\sum_{j\geq
    k}\|u^{j}\|_2^2+c(k^{7/2}\|u_0\|_2^3+\nu k^3\| u_0\|_2^2)$$
    and hence
$$\di\sum_{j\geq k}\| u^{j}(t)\|_2^2\leq c(k,\|u_0\|_2,\nu)\int_0^t\big(\sum_{j\geq
    k}\|u^{j}\|_2^2\big)\,d\tau+\sum_{j\geq
    k}\|u_0^{j}\|_2^2+e_1(k,\|u_0\|_2,\nu,T)$$
    for all $t\in (0,T)$.
Thus, by Gronwall's inequality
$$\sum_{j\geq k}\| u^{j}(t)\|_2^2\leq \Big(\sum_{j\geq
    k}\|u_0^{j}\|_2^2+e_1(k,\|u_0\|_2,\nu,T)\Big)\exp(c(k,\|u_0\|_2,\nu)T),
             \quad \forall t\in (0,T).$$
This inequality together with \eq{2.11} imply that $u^k\in
L^\infty(0,T;\dot{H}^{1/2}(\R^3))$, and hence $u$ can not blow up at
$t=T$. Thus the proof is complete.

\hfill \qed 
\section{Regularity in $L^\infty(0,T; B^{-1}_{\infty,\infty})$: Proof of the main result}

By a simple linear energy superposition, as seen in the proof of
Theorem \ref{T2.3}, one obtains the regularity criterion assuming
smallness of $L^\infty(B^{-1}_{\infty,\infty})$-norm for Leray Hopf
weak solutions. In this section we show that this smallness
condition can be removed, thus giving the proof of Theorem
\ref{T1.1}. The idea is to apply double superposition to the energy
estimates for high frequency parts. More precisely, we apply, first,
polynomial energy superposition with a large enough polynomial order
depending on the size of $L^\infty(B^{-1}_{\infty,\infty})$-norm and
then apply the linear superposition repeatedly to the resulting
differential inequality.
\par\bigskip\noindent
 {\bf Proof of Theorem \ref{T1.1}:}
 Let us assume without loss of generality that $T$ is
the first blow-up epoch for $u$. Let
 $$m:=\|u\|_{L^\infty(0,T;{B}^{-1}_{\infty,\infty})}.$$
Given any $k>0$, by \eq{2.19} and \eq{2.21} with $\si=-1$ we have
 \begin{equation}
 \label{E2.17}
\begin{array}{l}
  \frac{d}{2dt}\|u^{k}\|_2^2+\nu \|\na u^k\|_2^2
   \leq C_1mk^2\|u^{k/2}\|_2^2\ek
\qquad \leq C_1mk^2(\|u_{k/2,\tilde{m}
k}\|_2^2+\|u^{\tilde{m}k}\|_2^2)\ek \qquad \leq
C_1mk^2\|u_{k/2,\tilde{m}k}\|_2^2+C_1m \tilde{m}^{-2} \|\na
u^{\tilde{m}k}\|_2^2
\end{array}
 \end{equation}
for some absolute constant $C_1>0$ and any $\tilde{m}> 1$. Hence, if
$$\tilde{m}:=\sqrt{2C_1m\nu^{-1}},$$
 one has $C_1m \tilde{m}^{-2}\leq
\nu/2$ and
 \begin{equation}
 \label{E2.3}
 \begin{array}{rcl}
  \frac{d}{2dt}\|u^{k}\|_2^2+\frac{\nu}{2} \|\na u^k\|_2^2 &\leq &
  C_1mk^2\|u_{k/2,\tilde{m}k}\|_2^2\ek
 & \leq & \frac{C_2m}{2}\|\na u_{k/2,\tilde{m}k}\|_2^2.
\end{array}
 \end{equation}

From now on, let $k=j^s, j\in\N$, with positive integer $s$. Then,
\eq{2.3} implies that
 \begin{equation}
 \label{E2.4}
 \begin{array}{l}
  \frac{d}{dt}\|u^{j^s}\|_2^2+\nu\|\na u^{j^s}\|_2^2\leq
  C_2m\|\na u_{j^{s}/2,\tilde{m}j^{s}}\|_2^2.
\end{array}
 \end{equation}
Adding up \eq{2.4} over $j\in \N$ with $j\geq l$ for fixed $l\in\N$,
we have
 \begin{equation}
 \label{E2.5}
  \begin{array}{l}
 \di\frac{d}{dt}\sum_{j\geq l}\|u^{j^s}\|_2^2+\nu \sum_{j\geq l}
 \|\na u^{j^s}\|_2^2  \leq  C_2m\di \sum_{j\geq l}\|\na
 u_{j^{s}/2,\tilde{m}j^{s}}\|_2^2.
   \end{array}
 \end{equation}
Note that
 \begin{equation}
 \label{E2.25}
\begin{array}{l} \di \sum_{j\geq l}\|u^{j^s}\|_2^2
=\di\sum_{j\geq l}(j-l+1)\|u_{j^s,
 {(j+1)}^s}\|_2^2,\\[3ex]
 \di\sum_{j\geq l}\|\na u^{j^s}\|_2^2
 =\sum_{j\geq l}(j-l+1)\|\na u_{j^s,{(j+1)}^s}\|_2^2.
 \end{array}
 \end{equation}
Moreover,
 \begin{equation}
 \label{E2.26}
\begin{array}{l}
 \di\sum_{j\geq l}\|\na  u_{j^{s}/2,\tilde{m}j^{s}}\|_2^2
  =\di\sum_{j=l}^{[(2\tilde{m})^{1/s}l]-1}(j-l+1)\|\na
u_{j^{s}/2,(j+1)^{s}/2}\|_2^2\\[2ex]
\qquad +([(2\tilde{m})^{1/s}l]-l+1)\|\na
u_{[(2\tilde{m})^{1/s}l]^{s}/2,\tilde{m}l^{s}}\|_2^2\\[2ex]
\qquad +([[(2\tilde{m})^{1/s}l](1-(2\tilde{m})^{-1/s})]+1)\|\na
u_{\tilde{m}l^{s}, ([(2\tilde{m})^{1/s}l]+1)^{s}/2}\|_2^2\\[2ex]
 \qquad +\di\sum_{j\geq[(2\tilde{m})^{1/s}l]+1}([(1-(2\tilde{m})^{-1/s})j]+1)\|\na
u_{j^{s}/2,(j+1)^{s}/2}\|_2^2\\[4ex]
 \leq \|\na u^{l^s/2}\|_2^2+((2\tilde{m})^{1/s}-1)l\|\na
u_{l^{s}/2,\tilde{m}l^{s}}\|_2^2+(1-(2\tilde{m})^{-1/s})\times\\[2ex]
\quad \times\big([(2\tilde{m})^{1/s}l]\|\na u_{\tilde{m}l^{s},
([(2\tilde{m})^{1/s}l]+1)^{s}/2}\|_2^2
 + \di\sum_{j\geq [(2\tilde{m})^{1/s}l]+1}j\|\na
u_{j^{s}/2,(j+1)^{s}/2}\|_2^2\big),
  \end{array}
  \end{equation}
   where $[\cdot]$ denotes the integer part. Observe that
   $([(2\tilde{m})^{1/s}l]+1)^s/2\geq \tilde{m}l^s$ and that  $j^s\leq
(j+i(j))^s/2\leq (j+1)^s$ for some $i(j)\in\N$ iff $2^{1/s}j\leq
j+i(j)\leq 2^{1/s}(j+1)$. Hence
$$\begin{array}{l}
\di[(2\tilde{m})^{1/s}l]\|\na u_{\tilde{m}l^{s},
([(2\tilde{m})^{1/s}l]+1)^{s}/2}\|_2^2+
 \sum_{j\geq [(2\tilde{m})^{1/s}l]+1}j\|\na
 u_{j^{s}/2,(j+1)^{s}/2}\|_2^2\\[4ex]
 \qquad \leq
\di\sum_{j\geq l}(j+i(j))\|\na u_{j^{s},(j+1)^{s}}\|_2^2\leq
  \sum_{j\geq l}2^{1/s}(j+1)\|\na u_{j^{s},(j+1)^{s}}\|_2^2,
 \end{array}$$
where
 $$2^{1/s}(j+1)\leq (2\tilde{m})^{1/s}j,\quad \forall j\geq
 (\tilde{m}^{1/s}-1)^{-1}.$$
  Therefore, for $l\geq (\tilde{m}^{1/s}-1)^{-1}$ it follows from \eq{2.26} that
  \begin{equation}
 \label{E2.27}
 \begin{array}{l}
 \di\sum_{j\geq l}\|\na  u_{j^{s}/2,\tilde{m}j^{s}}\|_2^2
 \leq \|\na u^{l^s/2}\|_2^2
 +((2\tilde{m})^{1/s}-1)l \|\na u_{l^s/2,\tilde{m}l^s}\|_2^2\\[3ex]
\qquad   +((2\tilde{m})^{1/s}-1)\di\sum_{j\geq l}j\|\na u_{j^s,(j+1)^s}\|_2^2\\[3ex]
\quad \leq \|\na u^{l^s/2}\|_2^2
 +((2\tilde{m})^{1/s}-1)\big(l \|\na u_{l^s/2,\tilde{m}l^s}\|_2^2
   +\di\sum_{j\geq l}j\|\na u_{j^s,(j+1)^s}\|_2^2\big).
 \end{array}
 \end{equation}
Thus, it follows from \eq{2.5}, \eq{2.25} and \eq{2.27} that for
$l\geq (\tilde{m}^{1/s}-1)^{-1}$
\begin{equation}
 \label{E2.6}
  \begin{array}{l}
 \di\frac{d}{dt}\sum_{j\geq l}(j-l+1)\|u_{j^s, (j+1)^s}\|_2^2+\nu
  \sum_{j\geq l}(j-l+1)\|\na u_{j^s, (j+1)^s}\|_2^2 \\[3ex]
 \qquad \leq A(l)+a(s)\big(l \|\na u_{l^s/2,\tilde{m}l^s}\|_2^2
   +\di\sum_{j\geq l}j\|\na u_{j^s,(j+1)^s}\|_2^2\big),
  \end{array}
 \end{equation}
 where and in what follows
 $$A(l)\equiv C_2m \|\na u^{l^s/2}\|_2^2,\quad
   a(s)\equiv C_2m((2\tilde{m})^{1/s}-1).$$
Therefore, if
 \begin{equation}
 \label{E2.20}
a(s)<\frac{\nu}{4},
 \end{equation}
then it follows from \eq{2.6} that
$$\di\sum_{j\geq l}j\|u_{j^s, (j+1)^s}(t)\|_2^2+\int_0^T \sum_{j\geq
l}j\|\na u_{j^s, (j+1)^s}\|_2^2\leq c(m,\nu,l)\|u_0\|_2^2.
$$
Thus, in view of
$$j\|u_{j^s, (j+1)^s}\|_2^2\sim \|u_{j^s, (j+1)^s}\|_{\dot{H}^{1/(2s)}}^2,$$
 we have
\begin{equation}
\label{E2.29}
\begin{array}{l}
u\in L^\infty(0,T; H^{\frac{1}{2s}}(\R^3))\cap  L^2(0,T;
H^{1+\frac{1}{2s}}(\R^3)),\ek
 \|u\|_{L^\infty(0,T;
H^{\frac{1}{2s}}(\R^3))}+\|u\|_{L^2(0,T;
H^{1+\frac{1}{2s}}(\R^3))}\leq C(m,\nu)\|u_0\|_2.
\end{array}
\end{equation}
\par\medskip
Let us proceed the argument further. In view of $j-l+1\geq
\frac{1}{2}(j+1)$ for $j\geq 2l-1$, it follows from \eq{2.6} that
$$ \begin{array}{l}
 \di\frac{d}{dt}\sum_{j\geq l}(j-l+1)\|u_{j^s, (j+1)^s}\|_2^2+\frac{\nu}{2}
  \sum_{j\geq 2l-1}(j+1)\|\na u_{j^s, (j+1)^s}\|_2^2 \\[3ex]
 \quad \leq A(l)+a(s)(l \|\na u_{l^s/2,\tilde{m}l^s}\|_2^2
   +\di\sum_{j\geq l}j\|\na u_{j^s,(j+1)^s}\|_2^2)
  \end{array}$$
and thus, under the condition \eq{2.20},
 \begin{equation}
 \label{E2.8}
\begin{array}{l}
 \di\frac{d}{dt}\sum_{j\geq l}(j-l+1)\|u_{j^s, (j+1)^s}\|_2^2+\frac{\nu}{4}
  \sum_{j\geq 2l-1}(j+1)\|\na u_{j^s, (j+1)^s}\|_2^2 \\[3ex]
\quad \leq A(l)+a(s)(l \|\na u_{l^s/2,\tilde{m}l^s}\|_2^2
   +\di\sum_{j=l}^{2l-2}j\|\na u_{j^s,(j+1)^s}\|_2^2).
  \end{array}
 \end{equation}

Let us multiply \eq{2.8} with $l^{i-1}$ for $i\in\N$ with $i\leq
2s-1$ and add up the result over $l\in \N$ with $l\geq l_1\geq
(\tilde{m}^{1/s}-1)^{-1}$, we have
\begin{equation}
 \label{E2.7}
  \begin{array}{l}
 \di\frac{d}{dt}\sum_{l\geq l_1}\sum_{j\geq l}(j-l+1)l^{i-1}\|u_{j^s, (j+1)^s}\|_2^2
 +\frac{\nu}{4}
  \sum_{l\geq l_1}\sum_{j\geq 2l-1}(j+1)l^{i-1}\|\na u_{j^s, (j+1)^s}\|_2^2 \\[3ex]
\quad\di\leq \sum_{l\geq l_1}l^{i-1}A(l)+a(s)\Big(\sum_{l\geq
l_1}l^i \|\na u_{l^s/2,\tilde{m}l^s}\|_2^2
   +\di\sum_{l\geq l_1}\sum_{j=l}^{2l-2}jl^{i-1}\|\na u_{j^s,(j+1)^s}\|_2^2\Big).
  \end{array}
 \end{equation}

In the below, $P_{i}^{(k)}(j)$ denotes a polynomial in $j$ up to
$i$-th order, the coefficient of which may depend on $i$. In
\eq{2.7}, we have
 \begin{equation}
 \label{E2.39}
\begin{array}{l}
 \di \sum_{l\geq l_1}\sum_{j\geq l}(j-l+1)l^{i-1}\|u_{j^s,(j+1)^s}\|_2^2 \\[3ex]
 \quad  =\di\sum_{j\geq l_1}
  (l_1^{i-1}+\cdots +j^{i-1})(j+1)\| u_{j^s,(j+1)^s}\|_2^2\\[2ex]
\qquad -\di\sum_{j\geq l_1}
  (l_1^{i}+\cdots +j^{i})\|u_{j^s,(j+1)^s}\|_2^2\\[3ex]
 \quad  \geq \di\sum_{j\geq l_1}
  i^{-1}(j^i-l_1^{i})(j+1)\| u_{j^s,(j+1)^s}\|_2^2\\[2ex]
\qquad -\di\sum_{j\geq l_1}
 (i+1)^{-1}  ((j+1)^{i+1}-(l_1-1)^{i+1})\|u_{j^s,(j+1)^s}\|_2^2\\[3ex]
\quad  = (i^2+i)^{-1}\di\sum_{j\geq l_1}  j^{i+1}\|
u_{j^s,(j+1)^s}\|_2^2+\di\sum_{j\geq l_1}
P_{i}^{(1)}(j)\|u_{j^s,(j+1)^s}\|_2^2
 \end{array}
 \end{equation}
 and
\begin{equation}
 \label{E2.39n}
\begin{array}{l}
 \di \sum_{l\geq l_1}\sum_{j\geq l}(j-l+1)l^{i-1}\|u_{j^s,(j+1)^s}\|_2^2 \\[3ex]
\quad  \leq i^{-1}\di\sum_{j\geq
l_1}\big((j+1)^{i+1}-(l_1-1)^{i+1}\big)\|u_{j^s,(j+1)^s}\|_2^2\\[2ex]
 \qquad -(i+1)^{-1}\di\sum_{j\geq
 l_1}\big(j^{i+1}-l_1^{i+1}\big)\|u_{j^s,(j+1)^s}\|_2^2\\[3ex]
\quad =
 (i^2+i)^{-1}\di\sum_{j\geq l_1}  j^{i+1}\|
u_{j^s,(j+1)^s}\|_2^2+\di\sum_{j\geq l_1}
P_{i}^{(2)}(j)\|u_{j^s,(j+1)^s}\|_2^2.
 \end{array}
 \end{equation}
Here  we used that
$$\begin{array}{c}
l^{i-1}+(l+1)^{i-1}+\cdots+j^{i-1}=\sum_{p=1}^j
p^{i-1}-\sum_{p=1}^{l-1} p^{i-1}
  \end{array}$$
and
 \begin{equation}
  \label{E2.40}
  i^{-1}j^{i}\leq \sum_{p=1}^j p^{i-1}\leq i^{-1}(j+1)^{i},\quad
  \forall  i,j\in\N.
 \end{equation}
Similarly, we have
 \begin{equation}
 \label{E2.43}
\begin{array}{l}
\di\sum_{l\geq l_1}\sum_{j\geq 2l-1}(j+1)l^{i-1}\|\na u_{j^s,
(j+1)^s}\|_2^2\\[3ex]
 \quad  =\di\sum_{j\geq 2l_1-1}
  \big((2l_1-1)^{i-1}+\cdots +j^{i-1}\big)(j+1)\|\na u_{j^s,(j+1)^s}\|_2^2\\[3ex]
 \quad  \geq i^{-1}\di\sum_{j\geq 2l_1-1}  j^{i+1}\|\na u_{j^s,(j+1)^s}\|_2^2
  +\di\sum_{j\geq 2l_1-1} P_{i}^{(3)}(j)\|\na u_{j^s,(j+1)^s}\|_2^2.
 \end{array}
 \end{equation}
 It is clear that
\begin{equation}
\label{E2.13}
  \di\sum_{l\geq l_1}l^{i-1}A(l)\equiv C_2m \sum_{l\geq l_1}l^{i-1}\|\na u^{l^s/2}\|_2^2
=\sum_{j\geq l_1}P_{i}^{(4)}(j)\|\na u_{j^s/2,(j+1)^s/2}\|_2^2.
\end{equation}
On the other hand, using \eq{2.40}, we have
 \begin{equation}
 \label{E2.41}
\begin{array}{l}
 \di\sum_{l\geq l_1}l^i \|\na u_{l^s/2,\tilde{m}l^s}\|_2^2  \leq
  \di\sum_{j= l_1}^{[(2\tilde{m})^{1/s}l_1]}(l_1^{i}+\cdots+j^{i})
        \|\na  u_{j^{s}/2,(j+1)^{s}/2}\|_2^2\\[2ex]
 \qquad +\di \sum_{j\geq[(2\tilde{m})^{1/s}l_1]+1}
    ([(2\tilde{m})^{-1/s}j]^{i}+\cdots+ j^{i})\|\na  u_{j^{s}/2,(j+1)^{s}/2}\|_2^2\\[4ex]
\quad \leq \di (i+1)^{-1}\Big(\sum_{j=l_1}^{[(2\tilde{m})^{1/s}l_1]}
\big((j+1)^{i+1}-(l_1-1)^{i+1}\big)\|\na  u_{j^{s}/2,(j+1)^{s}/2}\|_2^2\\[2ex]
 \qquad +\di \sum_{j\geq[(2\tilde{m})^{1/s}l_1]+1}
    \big((j+1)^{i+1}-[(2\tilde{m})^{-1/s}j]^{i+1}\big)\|\na  u_{j^{s}/2,(j+1)^{s}/2}\|_2^2\Big)
      \\[4ex]
 \quad\leq \di  c_1(\tilde{m},l_1,i)\|\na u^{l_1^s/2}\|_2^2
  +\sum_{j\geq [(2\tilde{m})^{1/s}l_1]+1}P_i^{(5)}(j)\|\na u_{j^{s}/2,(j+1)^{s}/2}\|_2^2
  \\[2ex]
  \qquad +\di (i+1)^{-1}\big(1-(2\tilde{m})^{-(i+1)/s}\big)
                           \sum_{j\geq [(2\tilde{m})^{1/s}l_1]+1}
   j^{i+1}\|\na u_{j^{s}/2,(j+1)^{s}/2}\|_2^2.
  \end{array}
  \end{equation}
 Concerning the third term on the right-hand side of \eq{2.41}, we have
 \begin{equation}
 \label{E2.42}
\begin{array}{l}
 \di (i+1)^{-1}(1-(2\tilde{m})^{-(i+1)/s})
         \di\sum_{j\geq [(2\tilde{m})^{1/s}l_1]+1}j^{i+1}\|\na
 u_{j^{s}/2,(j+1)^{s}/2}\|_2^2\\[3ex]
 \quad \leq
  \di(i+1)^{-1}(1-(2\tilde{m})^{-(i+1)/s})
   \sum_{j\geq l_1}2^{(i+1)/s}(j+1)^{i+1}\|\na
   u_{j^{s},(j+1)^{s}}\|_2^2 \\[3ex]
\quad \leq 4(i+1)^{-1}\di\sum_{j\geq l_1}j^{i+1}\|\na
u_{j^{s},(j+1)^{s}}\|_2^2
    +\di\sum_{j\geq l_1}P_i^{(6)}(j)\|\na u_{j^{s},(j+1)^{s}}\|_2^2
  \end{array}
  \end{equation}
in view of $i\leq 2s-1$. Finally,
 \begin{equation}
 \label{E2.45}
 \begin{array}{l}
 \di\sum_{l\geq l_1}\sum_{j=l}^{2l-2}jl^{i-1}\|\na
 u_{j^s,(j+1)^s}\|_2^2\leq \sum_{j= l_1}^{2l_1-2}(l_1^{i-1}+\cdots+j^{i-1})j\|\na
 u_{j^s,(j+1)^s}\|_2^2\\[3ex]
 \qquad +\di\sum_{j\geq 2l_1-1}\big(([\frac{j}{2}]+1)^{i-1}+\cdots+j^{i-1}\big)j\|\na
 u_{j^s,(j+1)^s}\|_2^2\\[3ex]
 \quad \leq c_2(l_1,i)\|\na u^{l_1^s}\|_2^2
  +i^{-1}\di\sum_{j\geq 2l_1-1}\big((j+1)^i-[\frac{j}{2}]^{i}\big)
                                      j\|\na u_{j^s,(j+1)^s}\|_2^2\\[3ex]
  \quad \leq c_2(l_1,i)\|\na u^{l_1^s}\|_2^2
  +i^{-1}\di\sum_{j\geq 2l_1-1}(j^{i+1}+P_i^{(7)}(j))\|\na
  u_{j^s,(j+1)^s}\|_2^2.
 \end{array}
 \end{equation}

Integrating \eq{2.7} from $0$ to $t\in (0,T)$, we have
\begin{equation}
 \label{E2.44}
  \begin{array}{l}
 \di\sum_{l\geq l_1}\sum_{j\geq l}(j-l+1)l^{i-1}\|u_{j^s,
 (j+1)^s}(t)\|_2^2\\[2ex]
 +\di\frac{\nu}{4}\int_0^t
  \sum_{l\geq l_1}\sum_{j\geq 2l-1}(j+1)l^{i-1}\|\na u_{j^s, (j+1)^s}(\tau)\|_2^2\,d\tau\\[3ex]
\quad\di\leq \int_0^t\Big(\sum_{l\geq
l_1}l^{i-1}A(l)+a(s)\big(\sum_{l\geq l_1}l^i \|\na
u_{l^s/2,\tilde{m}l^s}(\tau)\|_2^2\\[2ex]
 +\di\sum_{l\geq l_1}\sum_{j=l}^{2l-2}jl^{i-1}\|\na
u_{j^s,(j+1)^s}\|_2^2\big)\Big)\,d\tau+\di\sum_{l\geq
l_1}\sum_{j\geq l}(j-l+1)l^{i-1}\|[u(0)]_{j^s,
 (j+1)^s}\|_2^2.
  \end{array}
 \end{equation}
Then  it follows from \eq{2.39}, \eq{2.43}--\eq{2.44} that for all
$t\in (0,T)$
\begin{equation}
 \label{E2.46}
  \begin{array}{l}
(i^2+i)^{-1}\di\sum_{j\geq l_1}  j^{i+1}\| u_{j^s,(j+1)^s}(t)\|_2^2
 +\frac{i^{-1}\nu}{4}\int_0^t \di\sum_{j\geq 2l_1-1}  j^{i+1}\|\na
 u_{j^s,(j+1)^s}\|_2^2\\[3ex]
  \quad \leq (i^2+i)^{-1}\di\sum_{j\geq l_1}  j^{i+1}\|
[u_0]_{j^s,(j+1)^s}\|_2^2+\di\sum_{j\geq l_1}
P_{i}^{(2)}(j)\|[u_0]_{j^s,(j+1)^s}\|_2^2\\
\qquad -\di\sum_{j\geq l_1} P_{i}^{(1)}(j)\|u_{j^s,(j+1)^s}(t)\|_2^2
 -\int_0^t \sum_{j\geq 2l_1-1} P_{i}^{(3)}(j)\|\na
 u_{j^s,(j+1)^s}(\tau)\|_2^2\,d\tau\\
 \qquad +\di\int_0^t
   \Big(\sum_{j\geq l_1}P_{i}^{(4)}(j)\|\na u_{j^s/2,(j+1)^s/2}(\tau)\|_2^2\\
    \qquad +a(s)\big(
\di  c_1(\tilde{m},l_1,i)\|\na u^{l_1^s/2}(\tau)\|_2^2
  +\sum_{j\geq [(2\tilde{m})^{1/s}l_1]+1}P_i^{(5)}(j)\|\na u_{j^{s}/2,(j+1)^{s}/2}(\tau)\|_2^2
  \\    \qquad +4(i+1)^{-1}\di\sum_{j\geq l_1}j^{i+1}\|\na
u_{j^{s},(j+1)^{s}}(\tau)\|_2^2 +\sum_{j\geq l_1}P_i^{(6)}(j)\|\na
    u_{j^{s},(j+1)^{s}}(\tau)\|_2^2\\
 \qquad +c_2(l_1,i)\|\na u^{l_1^s}(\tau)\|_2^2
  +i^{-1}\di\sum_{j\geq 2l_1-1}(j^{i+1}+P_i^{(7)}(j))\|\na
  u_{j^s,(j+1)^s}(\tau)\|_2^2\big)\Big)\,d\tau.
   \end{array}
 \end{equation}
Simplifying this inequality, we have
\begin{equation}
 \label{E2.47n}
  \begin{array}{l}
\frac{1}{i+1}\di\sum_{j\geq l_1}  j^{i+1}\| u_{j^s,(j+1)^s}(t)\|_2^2
 +\frac{\nu}{4}\int_0^t \di\sum_{j\geq 2l_1-1}  j^{i+1}\|\na
 u_{j^s,(j+1)^s}\|_2^2\,d\tau\\[3ex]
 \leq \frac{1}{i+1}\di\sum_{j\geq l_1}  j^{i+1}\|
[u_0]_{j^s,(j+1)^s}\|_2^2+\di\sum_{j\geq l_1}
P_{i}^{(2)}(j)\|[u_0]_{j^s,(j+1)^s}\|_2^2\\[2ex]
  -\di\sum_{j\geq l_1}
P_{i}^{(1)}(j)\|u_{j^s,(j+1)^s}(t)\|_2^2+5a(s)\int_0^t
   \di\sum_{j\geq 2l_1-1}j^{i+1}\|\na
u_{j^{s},(j+1)^{s}}\|_2^2\,d\tau\\
\qquad +\di\int_0^t\Big(\sum_{j\geq
l_1}P_{i}^{(8)}(j)\|\na u_{j^s/2,(j+1)^s/2}\|_2^2\\
  \qquad  +(a(s)+1)\di\sum_{j\geq l_1}P_{i}^{(9)}(j)\|\na u_{j^s,(j+1)^s}\|_2^2
     +c_3(l_1,i,s)\|\na u^{l_1^s}\|_2^2\Big)\,d\tau
   \end{array}
 \end{equation}
 for all $t\in (0,T)$ with some polynomials $P_{i}^{(8)}$ and $P_{i}^{(9)}$
 in $j$ up to $i$-th order.
 Thus it follows from \eq{2.47n} that, if
\begin{equation}
 \label{E2.47}
 a(s)<\frac{\nu}{40},
\end{equation}
then
\begin{equation}
\label{E2.48}
\begin{array}{l}
 \di\sum_{j\geq l_1}  j^{i+1}\|
u_{j^s,(j+1)^s}(t)\|_2^2
 +\frac{\nu(i+1)}{8}\int_0^t \di\sum_{j\geq 2l_1-1}  j^{i+1}\|\na
 u_{j^s,(j+1)^s}\|_2^2\,d\tau\\[3ex]
 \leq \di\sum_{j\geq l_1}  j^{i+1}\|
[u_0]_{j^s,(j+1)^s}\|_2^2+(i+1)\di\sum_{j\geq l_1}
P_{i}^{(2)}(j)\|[u_0]_{j^s,(j+1)^s}\|_2^2\\[2ex]
 \quad -(i+1)\di\sum_{j\geq l_1} P_{i}^{(1)}(j)\|u_{j^s,(j+1)^s}(t)\|_2^2
  +\di(i+1)\int_0^t\Big(\sum_{j\geq
l_1}P_{i}^{(8)}(j)\|\na u_{j^s/2,(j+1)^s/2}\|_2^2\\
  \quad  +(a(s)+1)\di\sum_{j\geq l_1}P_{i}^{(9)}(j)\|\na u_{j^s,(j+1)^s}\|_2^2
     +c_3(l_1,i,s)\|\na u^{l_1^s}\|_2^2\Big)\,d\tau
 \end{array}
\end{equation}
for all $t\in (0,T)$. In \eq{2.48}, note that
$$\begin{array}{l}
\di\sum_{j\geq l_1} j^{i+1}\|
[u_0]_{j^s,(j+1)^s}\|_2^2+(i+1)\di\sum_{j\geq l_1}
P_{i}^{(2)}(j)\|[u_0]_{j^s,(j+1)^s}\|_2^2\\[2ex]
\qquad\qquad \qquad\qquad  \leq
C(i)\|u_0\|_{H^{(i+1)/(2s)}(\R^3)}^2.
 \end{array}$$
Thus,  in view of
$$j^{i}\|u_{j^s, (j+1)^s}\|_2^2\sim \|u_{j^s, (j+1)^s}\|_{\dot{H}^{i/(2s)}}^2,$$
 for $i\leq 2s-1$ it follows from \eq{2.48}  that, under the
condition \eq{2.47},
 \begin{equation}
 \label{E2.49}
 u\in L^\infty(0,T; H^{(i+1)/(2s)}(\R^3))\cap  L^2(0,T;
H^{1+(i+1)/(2s)}(\R^3))
 \end{equation}
holds true provided that
 \begin{equation}
 \label{E2.50}
 u\in L^\infty(0,T; H^{i/(2s)}(\R^3))\cap  L^2(0,T;
H^{1+i/(2s)}(\R^3)).
 \end{equation}
By the way,  \eq{2.50} for $i=1$ holds under the condition \eq{2.47}
(see \eq{2.20} and \eq{2.29}). This means that \eq{2.49} holds for
$i=2s-1$ under the condition \eq{2.47}. We remark that one can find
such $s\in\N$ depending  explicitly only on $m$, $\nu$ and
independent of $i$ that \eq{2.47} is satisfied. That is, we obtain
 \begin{equation}
 \label{E2.51}
 u\in L^\infty(0,T; H^{1}(\R^3))\cap  L^2(0,T;H^{2}(\R^3)).
 \end{equation}
Consequently, $u$ can not blow up at $t=T$ and the proof is
complete.


\hfill\qed

\bigskip
\noindent {\bf Fund:} This work was supported by no fund
organization.

\bigskip\noindent
\textbf{Declaration of interest:} The author makes sure that there
is no actual or potential conflict with other people or organization
concerning the submission or publication of the manuscript.

\bigskip\noindent
\textbf{Data Availability:} No data were used to support this study.

\end{document}